\documentclass[twocolumn]{autart}    % Enable this line and disable the 
\usepackage{wrapfig}

\usepackage[utf8]{inputenc}   % codificação UTF-8
\usepackage[T1]{fontenc}      % suporte a caracteres acentuados

\usepackage{amsmath,amssymb,amsfonts,mathtools}
\usepackage{calrsfs,mathrsfs}
\usepackage{graphicx}
\usepackage{xcolor}
\usepackage{float}
\usepackage{multirow}
\usepackage{psfrag}
\usepackage{blindtext}

\usepackage{subcaption}
\newtheorem{assumption}{Assumption}

\newtheorem{theorem}{Theorem}
\newtheorem{lemma}{Lemma}

\begin{document}

\begin{frontmatter}

       %\title{Extremum Seeking Control for Actuation Dynamics Governed by a Wave PDE %with Distributed Effects\thanksref{footnoteinfo}}
       \title{Extremum Seeking Control for Wave-PDE Actuation with Distributed Effects}

\thanks[footnoteinfo]{This paper was not presented at any conference.\\
Corresponding author: Tiago Roux Oliveira.}

\author[PEL]{Elisio Juvenal Muchave} \ead{muchaveelisio@gmail.com},
\author[UERJ]{Pedro Henrique Silva Coutinho}\ead{phcoutinho@eng.uerj.br},\linebreak
\author[UERJ]{Tiago Roux Oliveira}\ead{tiagoroux@uerj.br},
\author[UCSD]{Miroslav Krsti\'{c}}\ead{mkrstic@ucsd.edu}

\address[PEL]{Graduate Program in Electronics Engineering, Rio de Janeiro State University (UERJ), Rio de Janeiro -- RJ, Brazil}
\address[UERJ]{Dept. of Electronics and Telecommunication Eng., Rio de Janeiro State University (UERJ), Rio de Janeiro -- RJ, Brazil}
\address[UCSD]{Dept. of Mechanical and Aerospace Engineering, University of California -- San Diego (UCSD), La Jolla -- CA, USA}

\begin{keyword}                           % Five to ten keywords,  
Extremum seeking; adaptive control; optimization; backstepping in infinite dimensions; partial differential equations; distributed wave models.
\end{keyword}                             % keyword list or with the 
                                          % help of the Automatica 
                                          % keyword wizard

\begin{abstract}                          % Abstract of not more than 200 words.
This paper deals with the gradient-based extremum seeking control (ESC) with actuation dynamics governed by distributed wave partial differential equations (PDEs). To achieve the control objective of real-time optimization for this class of infinite-dimensional systems, we first solve the trajectory generation problem to re-design the additive perturbation signal of the ESC system. Then, we develop a boundary control law through the backstepping method to compensate for the wave PDE with distributed effects, which ensures the exponential stability of the average closed-loop system by means of a Lyapunov-based analysis. At last, by employing the averaging theory for infinite-dimensional systems, we prove that the closed-loop trajectories converge to a small neighborhood surrounding the optimal point. Numerical simulations are presented to illustrate the effectiveness of the proposed method.
\end{abstract}

\end{frontmatter}

\section{Introduction}

\sloppy

Extremum seeking control (ESC) was introduced in~\cite{leblanc1922} in the context of maximizing power transfer in electrical systems. With the formulation proposed by \cite{Krsti2000StabilityOE}, ESC acquired a solid theoretical foundation and has been established as an effective real-time optimization strategy. ESC implements an iterative-adaptive approach to optimize nonlinear objective functions. Using the so-called gradient-based algorithm~\cite{Ariyur2003RealTimeOB}, 
%\cite{oliveira2021extremum,ghaffari2012multivariable,Manzie2009ExtremumSW},
this is performed by introducing small sinusoidal perturbations so that the output is gradually moved towards the optimal point, which can be a maximum or minimum of this objective function. ESC can thus be employed without requiring full knowledge of the map to be optimized~\cite{scheinker2024100}.

ESC has been employed in different theoretical frameworks, such as stochastic formulation~\cite{Manzie2009ExtremumSW}, multiparameter Newton-based generalization~\cite{Ghaffari2011MultivariableNE}, Lie-bracket extremum seeking~\cite{abdelgalil2024initialization}, unbiased extremum seeking~\cite{jbara2025constructive}, and event-triggered control implementations~\cite{rodrigues2025event}.
Moreover, ESC has also been considered in several applications, 
such as anti-lock brake systems, autonomous vehicles and mobile robots, renewable energy systems  
as well as bio-processes and refrigeration system optimization, see~\cite{bastin2009extremum,krstic2014extremum,matveev2015extremum,tan2010extremum} and references therein. However, this literature is limited to systems with finite-dimensional model representation.

Just recently, the ESC framework has been extended to deal with infinite-dimensional actuation dynamics governed by partial differential equations (PDEs)~\cite{Oliveira_Krstic_book_2022,oliveira2023extremum}. For addressing the infinite-dimensional nature of the actuation dynamics, it is necessary to properly design the additive probing signal of the ESC by solving the so-called trajectory generation problem \cite[Chapter-12]{Krsti2008BoundaryCO}, as well as develop a boundary controller to compensate for the propagation effect of the input signal through the PDE domain. Depending on the considered class of PDEs, this two-step recipe becomes difficult to solve. In~\cite{oliveira2016extremum}, the ESC for static maps with delays has been addressed by using a transport PDE for delay representation. This approach was later expanded for stochastic ESC with delayed dynamic maps~\cite{ruvsiti2018stochastic} and for non-constant delays~\cite{santos2020gradient}. Reference \cite{feiling2018gradient} has addressed the ESC with actuation dynamics governed by diffusion PDEs. In~\cite{oliveira2020multivariable}, the multivariable ESC was addressed, considering that each input is governed by a different diffusion PDE. This work also provided a generalization for different classes
of parabolic (reaction–advection–diffusion) PDEs, wave equations,
and first-order hyperbolic (transport-dominated) PDEs. The traffic congestion control with a downstream bottleneck was addressed by~\cite{yu2021extremum} using an ESC strategy along with the traffic density on the freeway segment described with a Lighthill–Whitham–Richards (LWR) macroscopic PDE. A complete guide for the ESC with wave PDE compensation has been discussed in~\cite{oliveira2021extremum}, with the extension for PDE--PDE cascades in~\cite{oliveira2021extremumPDEPDE}. In~\cite{silva2024extremum}, the ESC was addressed for a class of wave PDEs with {K}elvin-{V}oigt damping, and the publication \cite{yilmaz2024prescribed} has considered the prescribed-time ESC convergence for PDEs using chirpy probing. However, all the aforementioned works deal with PDEs without a distributed effect in the actuation dynamics.

In this context, the sole ESC publications considering distributed effects are \cite{Tsubakino_2023} and \cite{coutinho2025extremum}. While the authors in~\cite{Tsubakino_2023} have handled distributed delays (or transport PDEs), diffusion PDEs with distributed effects have been introduced in~\cite{coutinho2025extremum}. 
%In particular, the controller proposed in~\cite{bekiaris2011compensating} was adopted to %compensate for the distributed effects.
However, ESC with actuation dynamics governed by a wave PDE with distributed effects remains as an open problem, which motivates the studies to be presented here.  
Notably, this paper has been the first effort to pursue an extension
of ESC from the wave PDE to the class of \textit{distributed} wave PDEs. The complete design and analysis is carried out considering the scenario of Neumann-boundary actuation and Dirichlet-boundary sensing. To this end, we first provide a solution to the trajectory generation problem to obtain a suitable additive probing signal. This result is novel relative not only to \cite{Tsubakino_2023} and \cite{coutinho2025extremum} but also to other related work in the
literature. Second, inspired by the backstepping-like method proposed in \cite{Bekiaris_SCL}, the second contribution of this paper is to develop the compensation controller for the distributed wave PDE to ensure the asymptotic stability of the closed-loop average system by means of a Lyapunov-based stability analysis. Finally, it is shown that the trajectories of optimizing map converge to a neighborhood of the unknown extremum point, by employing the averaging theory for infinite-dimensional systems \cite{Hale1990AveragingII}. %~(see~Appendix~\ref{appendix:a}).

This paper is organized as follows. Section \ref{sec:formulation} presents the formulation of the ESC problem considering the actuator dynamics described by a wave PDE with distributed effects, including the definition of the demodulation and additive perturbation signal generation, as well as the steps to obtain the estimation error dynamics. Section~\ref{sec:main-result} describes the boundary control design via backstepping and the corresponding closed-loop averaged system. The stability and convergence analysis are presented in Section~\ref{sec:stability}. Section~\ref{sec:results} shows the numerical results that validate the proposed method, and Section~\ref{sec:conclusion} provides the concluding remarks with some directions for future work.

\textbf{Preliminaries:}~ 
The partial derivatives of a function $u(x,t)$ are denoted as 
$\partial_{x} u(x,t) = \partial u(x,t)/\partial x$, 
$\partial_{t} u(x,t) = \partial u(x,t)/\partial t$. 
% Whenever convenient, we use the compact form $u_{x}(x,t)$ e $u_{t}(x,t)$, or simply $u_{x}$ and $u_{t}$, respectively. 
The subscript ``$\mathrm{av}$'' denotes the average of a periodic variable with period $\Pi$. 
The Euclidean norm of a finite-dimensional state variable $\vartheta(t)$ is denoted as $|\vartheta(t)|$. 
The spatial norm in ${L}_{2}[0,D]$ of the PDE state $u(x, t)$ is denoted as 
$\|{u(t)}\|^{2}_{{L}_{2}([0,D])} = \int_{0}^{D} u^2(x,t) \,dx$. The index ${L}_{2}([0,D])$ is omitted for brevity, thus $\|{\cdot}\| = \|{\cdot}\|_{{L}_{2}([0,D])}$, if it is not specified. 
According to \cite{KH:02}, a vector-valued  function $f(t,\varepsilon) \in \mathbb{R}^{n}$ 
is said of order $\mathcal{O}(\varepsilon)$ over an interval $[t_{1},t_{2}]$, if 
$\exists k, \bar{\varepsilon}: |{f(t,\varepsilon)}| \leq k\varepsilon , \forall t \in [t_{1},t_{2}]$ and 
$\forall \varepsilon \in [0,\overline{\varepsilon}]$. In general, the precise estimates of the constants
$k$ and $\bar{\varepsilon}$ are not provided. Thus, we simply use $\mathcal{O}(\varepsilon)$, 
which is interpreted as an order of magnitude relation for a sufficiently small scalar $\varepsilon$.

\section{Extremum Seeking Control Formulation}
\label{sec:formulation}

\begin{figure}[!ht]
    \centering
    \includegraphics[width=\columnwidth]{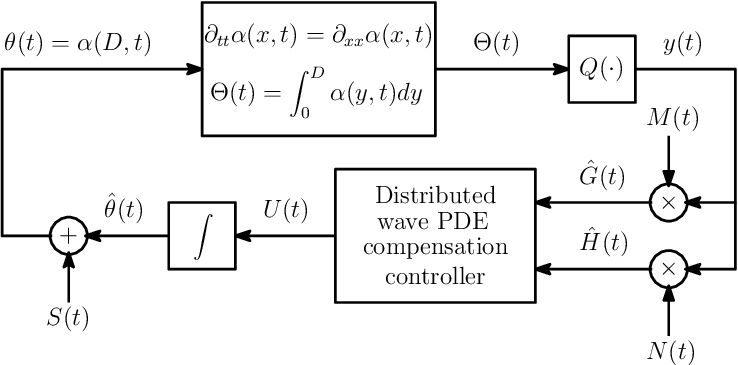}
    \caption{Extremum seeking control system with distributed wave PDEs.}
    \label{fig:esc_pde}
\end{figure}

Consider the gradient-based ESC system with actuation governed by a wave PDE with distributed effects shown in Figure~\ref{fig:esc_pde}. ESC aims to find and maintain the output $y \in \mathbb{R}$ of an unknown nonlinear static 
map $Q : \mathbb{R} \to \mathbb{R}$ within a small neighborhood of the unknown extremum point  $y^\ast$, corresponding to the optimizer $\Theta^\ast$, that is, $y^\ast = Q(\Theta^\ast)$.

\subsection{Nonlinear Maps with Wave-PDE Actuation and Distributed Effects}

The output of the map is given by
\begin{align}\label{eq:output}
    y(t) = Q(\Theta(t)),
\end{align}
where $\Theta(t)$ is the distributed input of the map, and the nonlinear map $Q$ is assumed to satisfy the following assumption.
\begin{assumption}\label{assump:quadratic_Q}
    The unknown static map is locally quadratic and described as:
    \begin{align}\label{eq:quadratic_map}
        Q(\Theta) = y^* + \frac{H}{2} \left( \Theta - \Theta^* \right)^2,
    \end{align}
    where $H  < 0$ is the unknown Hessian of the static map. 
\end{assumption}
Notice that Assumption~\ref{assump:quadratic_Q} is not restrictive since any continuously differentiable nonlinear map $Q$
can be locally approximated as a quadratic function around its extremum point. Without loss of generality, we consider $H < 0$ to address maximization problems. Minimization problems can be straightforwardly addressed by taking $H > 0$.
Thus, from \eqref{eq:output}--\eqref{eq:quadratic_map}, it follows that the output of the static map can be written as:
\begin{align}\label{eq:output_quadratic}
    y(t) = y^* + \frac{H}{2} \left( \Theta(t) - \Theta^* \right)^2.
\end{align}
We address the case in which the actuation dynamics is described by a distributed wave PDE given by 
\begin{align}
\Theta(t) &= \int_0^D \alpha(y,t) dy, \label{eq:6} \\ 
\partial_{tt} \alpha(x,t) &= \partial_{xx} \alpha(x,t), \label{eq:1} \\
\partial_{x} \alpha(0,t) &= 0, \label{eq:2b} \\
\alpha(D,t) &= \theta(t), \label{eq:2}
\end{align}
where $\alpha : [0,D] \times \mathbb{R}^+ \to \mathbb{R}$ is the state of the wave PDE, $D > 0$ corresponds to the upper bound of the spatial domain,
$\theta(t) \in \mathbb{R}$ is the actuation input and $\Theta(t)$ is the distributed actuation input of the nonlinear map $Q$ given in~\eqref{eq:quadratic_map}.

\subsection{Demodulation Signals}

Similar to \cite{Ghaffari2011MultivariableNE}, we consider the demodulation signal $M(t)$ to estimate the gradient of the static map
as follows:
\begin{equation}
\hat{G}(t) = M(t)y(t)\;\; \text{with} \;\; M(t) = \frac{2}{a}\sin{(\omega t)},
\label{eq:12}
\end{equation}
and the other demodulation signal $N(t)$ to estimate the Hessian of the static map as
\begin{equation}
\hat{H}(t) = N(t)y(t)\;\; \text{with} \;\; N(t) = -\dfrac{8}{a^{2}}\cos{(2\omega t)},
\label{eq:13}
\end{equation}
where $a>0$ and $\omega>0$ are parameters of amplitude and frequency, respectively.

\subsection{Trajectory Generation for the Motion Planning of the Additive Perturbation Signal}

Unlike the classical ESC, where the additive perturbation signal $S(t) = a \sin(\omega t)$ is directly applied to the input of the map, the presence of the distributed wave PDE significantly affects how the sinusoidal signal influences the input variable of the static map $Q(\cdot)$. This effect results from the distributed propagation of the signal over the spatial domain. To ensure that the sinusoidal perturbation is effectively present in the input of the map, it is necessary to design a boundary signal $S(t)$ such that, after being propagated by the distributed wave equation, the term $a \sin(\omega t)$ appears in the input of $Q(\cdot)$. In order to accomplish this goal, the following trajectory generation problem must be solved:
\begin{align}
    {\int_0^D \beta(y,t) dy }&= a\sin(\omega t),\label{eq:28} \\
    \partial_{tt}\beta(x,t) &= \partial_{xx}\beta(x,t),\label{eq:26}\\
    {\partial_x \beta(0,t)} &= 0,\label{27}\\
    \beta(D,t) &= S(t), \label{eq:25}
\end{align}
where $\beta : [0, D] \times \mathbb{R}^+ \rightarrow \mathbb{R}$, whose solution consists in determining $\beta(x, t)$ such that equations~\eqref{eq:28}--\eqref{eq:25} are satisfied. The following lemma states the analytical solution to the posed trajectory generation problem.
\begin{lemma}\label{lem:Sdesign}
    If the frequency $\omega > 0$ satisfies 
    \begin{align}\label{eq:probing_frequencies}
        \omega \neq \frac{k \pi}{D}, \quad k \in \mathbb{N},
    \end{align}
    then, the solution of the trajectory generation problem defined by equations~\eqref{eq:28}--\eqref{eq:25} is
    \begin{align}\label{eq:beta_solution}
        \beta(x,t) = \frac{a\omega}{\sin(\omega D)} \cos\left(\omega x\right) \sin(\omega t).
    \end{align}
    In particular, the additive perturbation signal $S(t)$ is  
    \begin{align}\label{eq:S_signal}
        S(t) = \beta(D,t) = \frac{a\omega}{\sin(\omega D)}  {\cos(\omega D)} \sin(\omega t).
    \end{align}
\end{lemma}
\begin{pf}
Consider a general solution in the following Taylor series form with time-varying coefficients:
\begin{align}\label{eq:candidate_beta}
    \beta(x,t) = \sum_{k=0}^{\infty} a_k(t) \frac{x^k}{k!}.
\end{align}
According to the boundary condition~\eqref{27}, we have that $a_1(t) = \partial_x\beta(0,t) = 0$.
Then, we substitute~\eqref{eq:candidate_beta} into~\eqref{eq:26}, which leads to
\begin{align}
    \sum_{k=0}^{\infty} \ddot{a}_k(t) \frac{x^k}{k!} = \frac{\partial}{\partial x^2} \sum_{k=0}^{\infty} a_k(t) \frac{x^k}{k!} 
     = \sum_{k=0}^\infty a_{k+2}(t) \frac{x^k}{k!},
\end{align}
from which we obtain the relation 
\begin{align}\label{eq:recurrency}
a_{k+2}(t) = \ddot{a}_k(t), 
\end{align}
which means that
$a_{2k+1}(t) = 0$, for $k \in \mathbb{N}$, and the general solution~\eqref{eq:candidate_beta} reduces to
\begin{align}\label{eq:candidate_beta_2}
     \beta(x,t) = \sum_{k=0}^{\infty} a_{2k}(t) \frac{x^{2k}}{(2k)!}.
\end{align}
Suppose that $a_0(t) = A \sin(\omega t) = A \mathrm{Im}\{e^{j\omega t}\}$, where $A$ is a constant to be determined.
From~\eqref{eq:recurrency}, we obtain that 
\begin{align}
    a_{2k}(t) = A (j \omega)^{2k} \mathrm{Im}\{e^{j\omega t} \}.
\end{align}
Then,~\eqref{eq:candidate_beta_2} becomes
\begin{align}\label{eq:candidate_beta_2cu}
   \beta(x,t) &= A \mathrm{Im}\left\lbrace e^{j\omega t} \sum_{k=0}^\infty \frac{(j \omega x)^{2k}}{(2k)!}  \right\rbrace \nonumber \\
   &= A \mathrm{Im} \{ e^{j \omega t} \cosh(j\omega x) \} = A \cos(\omega x) \sin(\omega t).
\end{align}
In addition, from~\eqref{eq:28}, we have that 
\begin{align}
    \int_0^D \beta(y,t) dy &= \int_0^D A \cos(\omega y) \sin(\omega t) dy \nonumber \\ 
    &=  A \frac{\sin(\omega D)}{\omega} \sin(\omega t) =  a\sin(\omega t),
\end{align}
from which follows that
\begin{align}\label{eq:A_solution}
    A = \frac{a\omega}{\sin(\omega D)},
\end{align}
which is finite if $\omega > 0$ is selected according to~\eqref{eq:probing_frequencies}.
Then, substituting~\eqref{eq:A_solution} into~\eqref{eq:candidate_beta_2cu}, we finally obtain $\beta(x,t)$ as in~\eqref{eq:beta_solution}. Evaluating $S(t) = \beta(D,t)$, we obtain~\eqref{eq:S_signal}. This concludes the proof. \hfill $\square$ \end{pf}

\subsection{Estimation Errors and Error Dynamics}

Consider the following estimate variables
\begin{align}
\hat{\theta}(t) = \theta(t) - S(t),\;\;\; \hat{\Theta}(t) = \Theta(t) - a\sin{(\omega t)},
\label{eq:14}
\end{align}
and the following estimation errors
\begin{align}
\tilde{\theta}(t) = \hat{\theta}(t) - \Theta^{*},\;\;\; \vartheta(t) = \hat{\Theta}(t) - \Theta^{*}.
\label{eq:15}
\end{align}
\newpage
According to~\eqref{eq:14}--\eqref{eq:15}, it is possible to establish the following relation between the distributed error $\vartheta(t)$, the distributed input $\Theta(t)$, and the optimum $\Theta^*$:
\begin{align}
    \vartheta(t) + \textcolor{black}{a\sin(\omega t)}= \Theta(t) - \Theta^{\ast}.\label{eq:24}
\end{align}
Let $\bar{\alpha}:[0,D]\times\mathbb{R}_{+}\rightarrow\mathbb{R}$ be defined as $\bar{\alpha}(x,t) = \alpha(x,t) - \beta(x,t) $.
From equations~\eqref{eq:6}--\eqref{eq:2} and \eqref{eq:28}--\eqref{eq:25}, along with \eqref{eq:14}--\eqref{eq:15}, it follows that:
\begin{align}
& \vartheta(t) = \int_0^D \bar{\alpha}(y,t) dy - \Theta^*, \label{eq:16} \\
& \partial_{tt}\bar{\alpha}(x,t) = \partial_{xx}\bar{\alpha}(x,t), \label{eq:17} \\
& \partial_{x}\bar{\alpha}(0,t) = 0, \label{eq:18} \\
& \bar{\alpha}(D,t) = \hat{\theta}(t). \label{eq:19}
\end{align}
The distributed error dynamics is obtained by taking the time derivative of equations \eqref{eq:16}--\eqref{eq:19} and using $\dot{\hat\theta} = \dot{\tilde{\theta}} = U(t)$ and $u(x,t) \coloneqq \partial_t \bar{\alpha}(x,t)$, such that:
\begin{align}
&    \dot{\vartheta}(t) = \int_0^D u(y,t) dy, \label{eq:20}\\
&    \partial_{tt}u(x,t) = \partial_{xx}u(x,t),\label{eq:21}\\
&    \partial_{x}u(0,t) = 0,\label{eq:22}\\
&    u(D,t) = U(t).\label{eq:23}
\end{align}

\subsection{Control Objective}

The main goal pursued in this paper is to design a distributed wave PDE compensation control law $U(t) = \Phi(\vartheta(t),u(x,t),\partial_t u(x,t))$ such that the average error dynamics of \eqref{eq:20}--\eqref{eq:23} is exponentially stable. Then, using the Average Theorem for infinite-dimensional systems \cite{Hale1990AveragingII}, 
%(see Appendix~\ref{appendix:a}),
prove the convergence of $(\theta(t), \Theta(t), y(t))$ to a neighborhood of the extremum $(\Theta^*, \Theta^*, y^*)$.

\section{Boundary Control and Closed-Loop System}
\label{sec:main-result}

\subsection{Backstepping Design in Infinite Dimensions}

Consider the distributed error dynamics described by equations~\eqref{eq:20}--\eqref{eq:23}. The following lemma provides a control design procedure to ensure the exponential stability of the closed-loop system with the proposed control law.
\begin{lemma}
\label{Lema:2}
Consider the closed-loop system resulting from the interconnection of~\eqref{eq:20}--\eqref{eq:23} with the controller
\begin{align}
U(t) &= \overline{K} Z(t) - c_0\int_0^D \partial_t u(y,t) dy ,\label{eq:30}
\end{align}
where
\begin{align}
Z(t) &= \vartheta(t) - g'(D) c_0  \int_0^D u(y,t) dy  \nonumber \\
    &+ \int_0^D g(y) \partial_t u(y,t) dy, \label{eq:31}
\end{align}
and
\begin{align}
    g(y) = \frac{D^2 - y^2}{2}, \label{eq:32}
\end{align}
for some $c_0 > 0$ with $c_0 \neq 1$.
If $\overline{K}$ is selected as a negative scalar such that
\begin{align}
\lambda = g'(D) \overline{K} > 0\label{eq:33}
\end{align}
satisfies
\begin{align}
    \lambda \neq \frac{1}{2D}\ln\left(-\frac{1 + c_0}{1 - c_0}\right),\label{eq:lambda_condition}
\end{align}
then, for any initial condition \( u(\cdot,0) \!\in\! H^1(0,D) \) and \( \partial_t u(\cdot,0) \!\in\! L^2(0,D) \), the closed-loop system\linebreak has a unique solution \( (\vartheta(t), u(\cdot,t), \partial_t u(\cdot,t)) \!\in\! \textcolor{black}{C\left([0,\infty), \mathbb{R}^n \times H^1(0,D) \times L^2(0,D)\right)}\), which is exponentially stable, in the sense that there exist positive constants $\kappa$ and $\lambda$ such that
\begin{align}
\Omega(t) &\leq \kappa \, \Omega(0) \, e^{-\rho t}, \label{eq:34}
\end{align}
with
\begin{align}
\Omega(t) = |\vartheta(t)|^2 + \|\partial_x u(t)\|^2 + \|\partial_t u(t)\|^2. \label{eq:35}
\end{align}
Moreover, if the initial condition $(u(\cdot,0), \partial_t u(\cdot,0))$ \textcolor{black}{is compatible} with the controller~\eqref{eq:30} and belongs to $H^2(0,D) \times H^1(0,D)$, then $(\vartheta(t), u(\cdot,t), \partial_t u(\cdot,t)) \in \textcolor{black}{C([0,\infty],\mathbb{R}^n \times H^1(0,D) \times L^2(0,D))}$ is the classical solution of the closed-loop system. 
\end{lemma}
\begin{pf}
Let us consider the invertible transformation applied to the finite-dimensional state \(\vartheta(t)\) given in~\eqref{eq:31} and the transformation for the infinite-dimensional actuator state \(u(x,t)\) defined as follows:
\begin{align}
    w(x,t) =\, & u(x,t) - \gamma(x)Z(t)  
    + c_0\int_0^x  \partial_t u(y,t) dy , \label{eq:36}
\end{align}
where the kernel \(\gamma(x)\) will be determined to transform the original system~\eqref{eq:20}--\eqref{eq:23}, together with the control law~\eqref{eq:30}, into the following target system:
\begin{align}
\dot{Z}(t) &= -\lambda Z(t), \label{eq:37} \\
\partial_{tt} w(x,t) &= \partial_{xx} w(x,t), \label{eq:38} \\
\partial_x w(0,t) &= c_0 \partial_t w(0,t), \label{eq:39} \\
 w(D,t) &= 0, \label{eq:40}
\end{align}
where \eqref{eq:33} is satisfied with $\overline{K} < 0$ since $g'(D) = -D$.

To verify this, we first differentiate \(Z(t)\) in~\eqref{eq:36}. Using the relation~\eqref{eq:21}, integration by parts, and the relations in~\eqref{eq:31} and \eqref{eq:32}, we obtain:
\begin{align}
\dot{Z}(t) &= \int_0^D u(y,t)  dy + g(D) \partial_x u(D,t) - g(0) \partial_x u(x,t) \nonumber \\
&- g'(D) u(D,t) + g'(0) u(0,t) + \int_0^D g''(y) u(y,t)  dy \nonumber \\
&- g'(D) c_0 \int_0^D \partial_{t} u(y,t) dy. \label{eq:41}
\end{align}
By noticing that the function \( g(\cdot) \) in~\eqref{eq:32}, used in the system transformation, is the solution to the following boundary value problem:
\begin{align}
g''(y) &= -1, \label{eq:42} \\
g'(0) &= 0, \label{eq:43} \\
g(D) &= 0, \label{eq:44}
\end{align}
we obtain
\begin{align}
\dot{Z}(t) = - g'(D) U(t) - g'(D) c_0 \int_0^D \partial_{t} u(y,t) \, dy. \label{eq:41b}
\end{align}
With the controller in~\eqref{eq:30} we obtain~\eqref{eq:37}.

In order to establish equations~\eqref{eq:38}--\eqref{eq:40}, we need to compute the derivatives of~\eqref{eq:36}. 
We initiate computing its first-order time derivative. Using \eqref{eq:21}, integration by parts, and~\eqref{eq:22}, we obtain 
\begin{align}
    \partial_{t} w(x,t) &= \partial_{t} u(x,t) + \lambda \gamma(x) Z(t) + c_0 \partial_x u(x,t). \label{eq:w_diff_t}
\end{align}
Then, the second-order time derivative of~\eqref{eq:36} is 
\begin{align}
\partial_{tt} w(x,t) = \partial_{tt} u(x,t) - \lambda^2 \gamma(x)  Z(t) + c_0 \partial_{xt}u(x,t). \label{eq:45}
\end{align}
Now we proceed to compute the spatial derivatives of~\eqref{eq:36}. The first spatial derivative is
\begin{align}
    \partial_{x} w(x,t) = \partial_{x} u(x,t) - \gamma'(x)Z(t) + c_0 \partial_{t}u(x,t) \label{eq:w_diff_x}
\end{align}
and the second spatial derivative of~\eqref{eq:36} is
\begin{align}
\partial_{xx} w(x,t) = \partial_{xx} u(x,t) - \gamma''(x) Z(t) + c_0 \partial_{tx} u(x,t). \label{eq:46}
\end{align}
From \eqref{eq:21}--\eqref{eq:23},~\eqref{eq:30}, \eqref{eq:36}, and \eqref{eq:w_diff_t}--\eqref{eq:46}, we have that
\begin{align}
w(0,t) &= u(0,t) - \gamma(0) Z(t), \label{eq:47} \\
\partial_{x} w(0,t) &= -\gamma'(0) Z(t) + c_0 \partial_{t} u(0,t), \label{eq:48} \\
\partial_{t} w(0,t) &= \partial_{t} u(0,t) + \lambda \gamma(0) Z(t), \label{eq:49} \\
w(D,t) &= \overline{K} Z(t) - \gamma(D) Z(t). \label{eq:50}
\end{align}
Using expressions~\eqref{eq:45}, \eqref{eq:46} and \eqref{eq:47}--\eqref{eq:50}, together with~\eqref{eq:21}--\eqref{eq:23}, we conclude that~\eqref{eq:38}--\eqref{eq:40} are satisfied if \( \gamma(x) \) is a solution of the following boundary value problem
\begin{align}
\gamma''(x) &= \lambda^2 \gamma(x), \label{eq:51}\\
\gamma'(0) &= - c_0 \lambda \gamma(0),\label{eq:5}\\
\gamma(D) &= \overline{K},\label{eq:53}
\end{align}
yielding to the following solution:
% \begin{align}
% \gamma(x) &= \alpha e^{A_{cl} x} + \beta e^{-A_{cl} x}\label{eq:54}
% \end{align}
\begin{align}
\gamma(x) &= \overline{K} \frac{e^{\lambda x} + r e^{-\lambda x}}{e^{\lambda D} + r e^{-\lambda D}}, \label{eq:54}
\end{align}
with 
\begin{align}
    r &= \frac{1 + c_0}{1 - c_0},
\end{align}
which exists for $c_0 \neq 1$ and $\lambda$ satisfying~\eqref{eq:lambda_condition}.

Consider now the Lyapunov functional candidate
\begin{align}
V(t) &= \frac{1}{2} Z^2(t) + E(t), \label{eq:62}
\end{align}
where the term $E(t)$, inspired by~\cite{smyshlyaev2009boundary}, is chosen as 
\begin{align}
E(t) &= \frac{1}{2} \Biggl(  \|\partial_x w(t)\|^2  + \|\partial_t w(t)\|^2 \Biggr) \nonumber \\ 
&+ \delta \int_0^D (y-2) \partial_y w(y,t) \partial_t w(y,t) dy,
\label{eq:64}
\end{align}
and $\delta$ is a sufficiently small positive scalar so that we ensure that $V(t)$ is positive definite. Differentiating~\eqref{eq:62} with respect to time, we obtain 
\begin{align}
\dot{V}(t) &= -\lambda Z^2(t) - \frac{\delta}{2}\left(\|\partial_x w(t)\|^2 + \|\partial_t w(t)\|^2 \right) \nonumber \\
& - [c_0 - \delta(1-c_0^2)]\partial_t w^2(0,t) - \frac{\delta}{2} \partial_x w^2(D,t) \nonumber \\
&\leq -\lambda Z^2(t) - \eta E(t) \nonumber \\
&\leq - \rho V(t), \label{eq:V1dot}
\end{align}
for some $\eta > 0$ and $\rho = \min\{\lambda,\eta\}$.
Then, using the comparison principle \cite{KH:02}, we get $V(t) \leq V(0) e^{-\rho t}$.

{\color{black}
To prove that \eqref{eq:34}--\eqref{eq:35} hold, it is sufficient to show that
\begin{align}
\underline{M} \Omega(t) \leq V(t) \leq \overline{M} \Omega(t),\label{eq:V_ineq}
\end{align}
for some positive constants \( \overline{M} \) and \( \underline{M} \).

% where $\Omega(t) = |\vartheta(t)|^2 + \|\partial_x u(t)\|^2 + \|\partial_t u(t)\|^2$. 

Using~\eqref{eq:31} and \eqref{eq:32} together with the Poincare, Young, and Cauchy-Schwarz's inequalities, we can conclude that there exists a positive constant $m$ such that
\begin{align}\label{eq:Z_bound}
    |Z(t)|^2 \!\leq\! m\left(|\vartheta(t)|^2 \!\!+\!\! \|\partial_x u(t)\|^2 \!\!+\!\! \|\partial_t u(t)\|^2\right).
\end{align}
Then, from \eqref{eq:62}, \eqref{eq:64}, and \eqref{eq:w_diff_t}, \eqref{eq:w_diff_x}, and \eqref{eq:Z_bound}, together with Young and Cauchy-Schwarz's inequalities, we obtain the upper bound in \eqref{eq:V_ineq}. The lower bound is obtained similarly by expressing $u(x,t)$ and $\vartheta(t)$ in terms of $w(x,t)$ (and its derivatives) and $Z(t)$.  The rest of the proof follows similar steps as~\cite{bekiaris2011compensating,smyshlyaev2009boundary}. This completes the proof.}
\hfill $\square$
% {\color{blue}
% Using transformations \eqref{eq:31} and \eqref{eq:36}, we can express $V(t)$ in terms of the original variables.
% From the definition of $w(x, t)$, we have:
% \begin{align} \label{ano2026}
%     \|w(t)\|^2 \!=\! \int_0^D \!\!\left( u(x, t) \!-\! \gamma(x)Z(t) \!+\! c_0 \int_0^x \!\!\partial_t u(y, t)  dy \right)^2\!dx.
% \end{align}
% %
% Applying the Cauchy--Schwarz and Young inequalities to (\ref{ano2026}), we get 
% \begin{align}
%     \|w(t)\|^2 \leq C_1\left( \|u(t)\|^2 + Z^2(t) + \|\partial_t u(t)\|^2 \right), 
% \end{align}
% %
% with $C_1 > 0$ dependent only on $c_0$, $D$, and $\sup_x |\gamma(x)|$. From expression \eqref{eq:31}, and employing again the Cauchy--Schwarz and Young inequalities, we also obtain 
% %
% \begin{align}
%     |Z(t)|^2 \!\leq\! 3\left(|\vartheta(t)|^2 \!\!+\!\! D^2 c_0^2 \|u(t)\|^2 \!\!+\!\! \|g\|_{L^2}^2 \|\partial_t u(t)\|^2\right).
% \end{align}
% Combining these estimates into $V(t) = \frac{1}{2}(Z^2 + \|w\|^2)$, it follows that
% \begin{align}
%     V(t) \leq \overline{M}\left(|\vartheta(t)|^2 + \|u(t)\|^2 + \|\partial_t u(t)\|^2\right),
% \end{align}
% for some $\overline{M} > 0$. Using once more \eqref{eq:31} and \eqref{eq:36} and the same inequalities, the following inverse estimate is obtained:
% \begin{align}
%     V(t) \geq \underline{M}\left(|\vartheta(t)|^2 + \|u(t)\|^2 + \|\partial_t u(t)\|^2\right)
% \end{align}
% %
% for some $\underline{M} > 0$ dependent on the same transformation constants. Thus, \eqref{eq:V_ineq} is satisfied and (\ref{eq:34}) can be obtained, with $\Omega(t)$ in (\ref{eq:35}).} 
% \hfill $\square$ 
\end{pf}

\subsection{Average Closed-Loop System}
\label{sec:Cav_design_ano_novo}

With the controller proposed in~\eqref{eq:30}, the distributed effect of the wave PDE dynamics can be effectively compensated. 
However, this controller cannot be directly implemented in real-time because it depends on $\vartheta(t)$, which is unavailable for measurement.
Hence, we employ similar arguments as in~\cite{Ghaffari2011MultivariableNE} to obtain an implementable controller for the ESC system with distributed wave PDE using averaging-based estimates of the gradient and Hessian. With this paradigm, we develop a real-time implementable version of the controller~\eqref{eq:30}.

First, by substituting~\eqref{eq:24} into~\eqref{eq:output_quadratic}, we obtain
\begin{align}\label{eq:71}
    y(t) = y^{\ast} + \frac{H}{2}(\vartheta(t)+a\sin{(\omega t)})^2.
\end{align}
Then, by substituting~\eqref{eq:71} into~\eqref{eq:12}, we have
\begin{align}
    \hat{G}(t) &= \frac{2}{a}\sin{(\omega t)}\left( y^{\ast} + \frac{H}{2}(\vartheta(t)+a\sin{(\omega t)})^2\right) \nonumber \\
        &= \frac{2}{a}y^{\ast}\sin{(\omega t)} + \frac{H}{a}\vartheta^2(t)\sin{(\omega t}) \nonumber \\
        &+ \frac{H}{a}\left(
        2a \vartheta(t)\sin^2{(\omega t)} + a^2 \sin^3{(\omega t)}\right)\label{eq:72}.
\end{align}
Finally, by substituting~\eqref{eq:71} into~\eqref{eq:13}, it results in
{
\begin{align}
    &\hat{H}(t) = -\frac{8}{a^2}\cos{(2\omega t)}\left(y^{\ast} + \frac{H}{2}(\vartheta(t)+a\sin{(\omega t)})^2\right) \nonumber \\
    &= -\frac{8}{a^2}y^{\ast}\cos{(2\omega t)}
    - \frac{4H}{a^2}\left[\vartheta^2(t)\cos{(2\omega t)}\right. \nonumber \\
    &+ \left.2a\vartheta(t)\sin{(\omega t)}\cos{(2\omega t)} + a^2\sin^2{(\omega t)}\cos{(2\omega t)} \right]\label{eq:73}.
\end{align}}
$\!\!$Consider the average of $\hat{G}(t)$ and $\hat{H}(t)$, over the time period of $2\pi/\omega$, that is, 
\begin{align}\label{eq:averaging}
    {\hat{G}_{\mathrm{av}}(t) = \frac{\omega}{2 \pi} \int_0^{\frac{2\pi}{\omega}} \hat{G}(t)dt,
    \; 
    \hat{H}_{\mathrm{av}}(t) = \frac{\omega}{2 \pi} \int_0^{\frac{2\pi}{\omega}} \hat{H}(t)dt}.
\end{align}
From~\eqref{eq:averaging}, we obtain the following average versions of the gradient and the Hessian estimates, given, respectively, by
\begin{align}\label{eq:75}
    \hat{G}_{\mathrm{av}}(t)  = H \vartheta_{\mathrm{av}}(t) \quad 
    \mathrm{and} \quad
    \hat{H}_{\mathrm{av}}(t) = H.
\end{align}
\noindent
Now, by substituting \eqref{eq:31} into \eqref{eq:30} and taking $\overline{K} = KH$, the control law becomes
\begin{align}\label{eq:88.}
U_{\mathrm{av}}(t) &= K \ \Big[\vartheta(t) - g'(D) c_0  \int_0^D u(y,t) dy \nonumber \\
&+ \int_0^D g(y) \partial_t u(y,t) dy \, dy\Big] 
 - c_0\int_0^D \partial_t u(y,t) dy,
\end{align}
with $K>0$. 

Assuming that the frequency $\omega$ is sufficiently large, we can consider $\vartheta(t)$ to be approximately constant over one averaging interval. 
Consequently, the average version of the controller~\eqref{eq:31} becomes:
\begin{align}
U_{\mathrm{av}}(t) &= K  \hat{G}_{\mathrm{av}}(t) \!+\! K\hat{H}_{\mathrm{av}}(t) \Big[ \!-\! g'(D) c_0  \int_0^D u_{\mathrm{av}}(y,t) dy  \nonumber \\
&+ \int_0^D g(y) \partial_t u_{\mathrm{av}}(y,t) dy\Big] - c_0\int_0^D \partial_t u_{\mathrm{av}}(y,t) dy.
\label{eq:75..}
\end{align}

Now, we can compute the average version of the distributed error dynamics~\eqref{eq:20}--\eqref{eq:23} as follows:
\begin{align}
 \dot{\vartheta}_{\mathrm{av}}(t) &= %kH\vartheta_{\mathrm{av}}t)+
 \int_0^{D}u_{\mathrm{av}}(y,t)dy, \label{eq:78}\\
    \partial_{tt} u_{\mathrm{av}}(x,t) &= \partial_{xx} u_{\mathrm{av}}(x,t),  \label{eq:79} \\
    \partial_x u_{\mathrm{av}}(0,t) &= 0, \label{eq:80} \\
    % \frac{d}{dt}u_{\mathrm{av}}(D,t) &= - cu_{\mathrm{av}}(D,t) + cU_{\mathrm{av}}(t), \label{eq:81}
    u_{\mathrm{av}}(D,t) &= U_{\mathrm{av}}(t), \label{eq:81}
\end{align}
where $\vartheta_{\text{av}}(t),\, U_{\text{av}}(t) \in \mathbb{R}$ and $x \in [0, D]$. Based on Lemma~\ref{Lema:2} and \textcolor{black}{following similar steps in \cite{K:2008_vascao}}, we conclude that the closed-loop averaged system governed by \eqref{eq:78}--\eqref{eq:81}  with the controller in \eqref{eq:75..}
% equation~\eqref{filteredcontrolreal2026}, or its equivalent average version~\eqref{eq:controlador_non_average_3} 
is exponentially stable.

\subsection{Implementable Control Law}
\label{sec:Cav_design}

Due to technical reasons in the application of the Average Theorem for infinite-dimensional systems \cite{Hale1990AveragingII} in the stability proof (\textcolor{black}{as well as inverse optimality purposes \cite{FOK:2023}}), we introduce a low-pass filter (with $c>0$) into the controller so that $U(t)$ can be treated as a state variable. 
% Finally, we obtain the averaged infinite-dimensional control law to compensate for the wave process.
% The controller dynamics are given by:
% \begin{align}
% U(t) = \frac{c}{s + c} U_{\text{av}}(t)
% \end{align}
% which can be rewritten in the Laplace domain as:
% \begin{align}
% (s + c) U(t) = c U_{\text{av}}(t)
% \end{align}
% Converting to the time domain, we obtain the differential equation:
% \begin{align}
% \frac{d}{dt} U(t) + c U(t) = c U_{\text{av}}(t)
% \end{align}
% Substituting the expression for $U_{\text{av}}(t)$, we arrive at the final control law:
Thus, inspired by the average version of the controller in~\eqref{eq:75..} 
using the fact that
\begin{align}
\int_0^D g(y) \partial_t u(y, t)  dy = - D \hat{\theta}(t) - \Theta(t) + a \sin(\omega t),
\end{align}
and the low-pass filter, we consider the following implementable controller:
\begin{align} \label{filteredcontrolreal2026}
&\frac{d}{dt} U(t) = - c U(t) + c\bigg\{ K \hat{G}(t)  \nonumber \\
&+ K \hat{H}(t) \bigg[ -g'(D) c_0 \int_0^D u(y,t)  dy - D \hat{\theta}(t) \nonumber \\
&- \Theta(t) + a \sin(\omega t) \bigg] - c_0 \int_0^D \partial_t u(y,t)  dy \bigg\}.
\end{align}
Therefore, the closed-loop system becomes
\begin{align}
&    \dot{\vartheta}(t) = \int_0^D u(y,t) dy, \label{eq:20_2}\\
&    \partial_{tt}u(x,t) = \partial_{xx}u(x,t),\label{eq:21_2}\\
&    \partial_{x}u(0,t) = 0,\label{eq:22_2}\\
&    \frac{d}{dt} u(D,t)=  - c u(D,t) + c\bigg\{ K \hat{G}(t)  \nonumber \\
&+ K \hat{H}(t) \bigg[ -g'(D) c_0 \int_0^D u(y,t)  dy - D \hat{\theta}(t) \nonumber \\
&- \Theta(t) + a \sin(\omega t) \bigg] - c_0 \int_0^D \partial_t u(y,t)  dy \bigg\}.\label{eq:23_2}
\end{align}

\section{Stability and Convergence Analysis}
\label{sec:stability}

The following theorem provides the local stability and convergence properties of the closed-loop system.

\begin{theorem} \label{teor}
Consider the system depicted in Figure~\ref{fig:esc_pde}, whose dynamics are described by a quadratic nonlinear map in Assumption~\ref{assump:quadratic_Q}, coupled with an actuation dynamics governed by the distributed wave equation defined in \eqref{eq:6}--\eqref{eq:2}. 
There exists a sufficiently large constant \( c > 0 \) such that, for some \(\bar{\omega}(c) > 0\), for every \(\omega > \bar{\omega}\) \textcolor{black}{satisfying~\eqref{eq:probing_frequencies}}, the closed-loop system described by
(\ref{eq:20_2})--(\ref{eq:23_2}), with control law (\ref{filteredcontrolreal2026}),  
%\eqref{eq:78}--\eqref{eq:81}
admits a unique periodic solution, locally exponentially stable in \(t\), with period \(\Pi := \frac{2\pi}{\omega}\), denoted by \(\vartheta^{\Pi}(t), u^{\Pi}(x,t)\), which satisfies the following inequality and ultimate bounds $\forall t \geq 0$:
\begin{align}
&\Big( \|\vartheta^{\Pi}(t)\|^2 + \|\partial_t u^{\Pi}(t)\|^2 
+ \|\partial_x u^{\Pi}(t)\|^2 \Big)^{1/2} 
\leq O\left( \frac{1}{\omega} \right), \label{vasco} \\
&\limsup_{t \to \infty} \left| \theta(t) - \Theta^* \right| = \mathcal{O}\left(\textcolor{black}{a\omega |\cot(\omega D)|} + \frac{1}{\omega} \right),\label{eq:output1} \\
&\limsup_{t \to \infty} \left| \Theta(t) - \Theta^* \right| = \mathcal{O}\left(a + \frac{1}{\omega} \right),\label{eq:output2} \\
&\limsup_{t \to \infty} \left| y(t) - y^* \right| = \mathcal{O}\left(a^2 + \frac{1}{\omega^2} \right) \label{eq:output3}.
\end{align}
\end{theorem}
\begin{pf}
From the exponential stability established in Lemma~\ref{Lema:2} and the invertibility of the transformation in~\eqref{eq:31}, it follows from (\ref{eq:34}) and (\ref{eq:35}) that
\begin{align}
|\vartheta_{\mathrm{av}}(t)|^2 
&+ \|\partial_t u_{\mathrm{av}}(t)\|^2 
+ \|\partial_x u_{\mathrm{av}}(t)\|^2 \notag \\
&\leq \textcolor{black}{\kappa} e^{- \textcolor{black}{\rho} t} \big( 
|\vartheta_{\mathrm{av}}(0)|^2 
+ \|\partial_t u_{\mathrm{av}}(0)\|^2 \notag \\
&\quad + \|\partial_x u_{\mathrm{av}}(0)\|^2 
\big), \quad \forall t > 0,
\label{eq:output4}
\end{align}
%with positive constants \( M, \overline{M} > 0 \).
after assuming $c \to +\infty$ in (\ref{filteredcontrolreal2026}) for simplicity, since (\ref{filteredcontrolreal2026}) results in (\ref{eq:30}). From~\eqref{eq:output4}, the exponential stability of the average system associated with the cascade of the ODE system (\ref{eq:78}) with the distributed wave PDE (\ref{eq:79})--(\ref{eq:81}) is guaranteed.  
%, where the energy decay occurs in the variables relevant to the hyperbolic context, %namely the temporal and spatial derivatives of the field \( u(x,t) \), as well as %the scalar variable \( \vartheta(t) \).
%
%From~\eqref{eq:output4}, the origin of the average closed-loop system, represented %in the original ODE PDE state variable \( u(x,t) \), is also exponentially stable.

Next, according to the averaging theory for infinite-dimensional systems~\cite{Hale1990AveragingII}, %(see Appendix~\ref{appendix:a}),
for \( \omega \) sufficiently large, the closed-loop system given by~\eqref{eq:20}–\eqref{eq:23}, with \( U(t) \) defined in~\eqref{filteredcontrolreal2026}, admits a unique exponentially stable periodic solution, around its equilibrium (origin), satisfying~\eqref{vasco}.

On the other hand, the asymptotic convergence to a neighborhood of the extremum point is established by taking the absolute value of both sides of (\ref{eq:24}),
resulting in: 
%
%taking the absolute value of the second expression in~\eqref{eq:14}, after %substituting \( \widehat{\Theta}(t) = \vartheta(t) + \Theta^* \), as defined %in~\eqref{eq:15}, resulting in:
\begin{align}
|\Theta(t) - \Theta^*| = |\vartheta(t) + a \sin(\omega t)|.
\label{eq:output5}
\end{align}
Considering~\eqref{eq:output5} and rewriting it by adding and subtracting the periodic solution \( \vartheta^{\Pi}(t) \), it follows that:
\begin{align}
|\Theta(t) - \Theta^*| = |\vartheta(t) - \vartheta^{\Pi}(t) + \vartheta^{\Pi}(t) + a \sin(\omega t)|.
\label{eq:output6}
\end{align}
By applying the Averaging Theorem~\cite{Hale1990AveragingII}, it can be concluded that 
\( \vartheta(t) - \vartheta^{\Pi}(t) \to 0 \) exponentially. 
Consequently, it follows that:
\begin{align}
\limsup_{t \to \infty} |\Theta(t) - \Theta^*| 
= \limsup_{t \to \infty} |\vartheta^{\Pi}(t) + a \sin(\omega t)|.
\label{eq:output7}
\end{align}
Finally, by using the relation in~\eqref{vasco}, the result presented in~\eqref{eq:output2} is obtained.
Provided that \( \theta(t) - \Theta^* = \tilde{\theta}(t) + S(t) \), as defined in~\eqref{eq:14} and \eqref{eq:15}, and recalling that \( S(t) \) in (\ref{eq:S_signal}) is of order \( \mathcal{O}(\textcolor{black}{a\omega |\cot(\omega D)|})\), 
with the following ultimate bound  
%
%as shown in~\eqref{eq:output2}, it follows that:
\begin{align}
\limsup_{t \to \infty} |\tilde{\theta}(t)| = \mathcal{O}\left(\frac{1}{\omega}\right), \label{eq:output8}
\end{align}
we obtain \eqref{eq:output1}.
To show the convergence of the output \( y(t) \), one can follow the same steps used for \( \Theta(t) \) by substituting~\eqref{eq:output2} into~\eqref{eq:output_quadratic}, so that:
\begin{align}
\limsup_{t \to \infty} |y(t) - y^*| = \limsup_{t \to \infty} |H \vartheta^{2}(t) + H a^{2} \sin^{2}(\omega t)|\label{eq:output9}.
\end{align}
Therefore, by rewriting~\eqref{eq:output9} in terms of \( \vartheta^{\Pi}(t) \), and again using (\ref{vasco}), one finally obtains~\eqref{eq:output3}. \hfill $\square$
\end{pf}

\section{Numerical Results}
\label{sec:results}

This section presents numerical simulations to illustrate the effectiveness of the proposed methodology.
Consider a static map~\eqref{eq:quadratic_map}
with Hessian $H\!=\!-2$, optimizer $\Theta^{\ast}\!=\!2$ and optimum output $y^{\ast}\!=\!5$. We consider the domain $D = 1$ for the wave PDE with distributed effects in~\eqref{eq:6}--\eqref{eq:2} describing the dynamics of the actuator. 
For the controller in~\eqref{filteredcontrolreal2026}, we consider the gain $K=0.1$ and the corner frequency $c = 10$ for the corresponding low-pass filter.
For the perturbation and demodulation signals, it is considered $a = 0.1$ and $\omega = 7.5~\mathrm{rad/s}$, which also satisfies the condition~\eqref{eq:probing_frequencies}. 
%For implementation purposes, similar to~\cite{Ghaffari2011MultivariableNE}, we also %consider a high-pass filter in the output with cut frequency~$\omega_h = %0.75$~rad/s, a first-order low-pass filter in the gradient estimation $\hat{G}(t)$, %and a second-order low-pass filter in the Hessian estimation $\hat{H}(t)$, both with %cut frequency $\omega_l = 0.75$~rad/s. 

\begin{figure}[!ht]
    \centering
    \begin{subfigure}[b]{\columnwidth}
        \centering
        \includegraphics[width=\columnwidth]{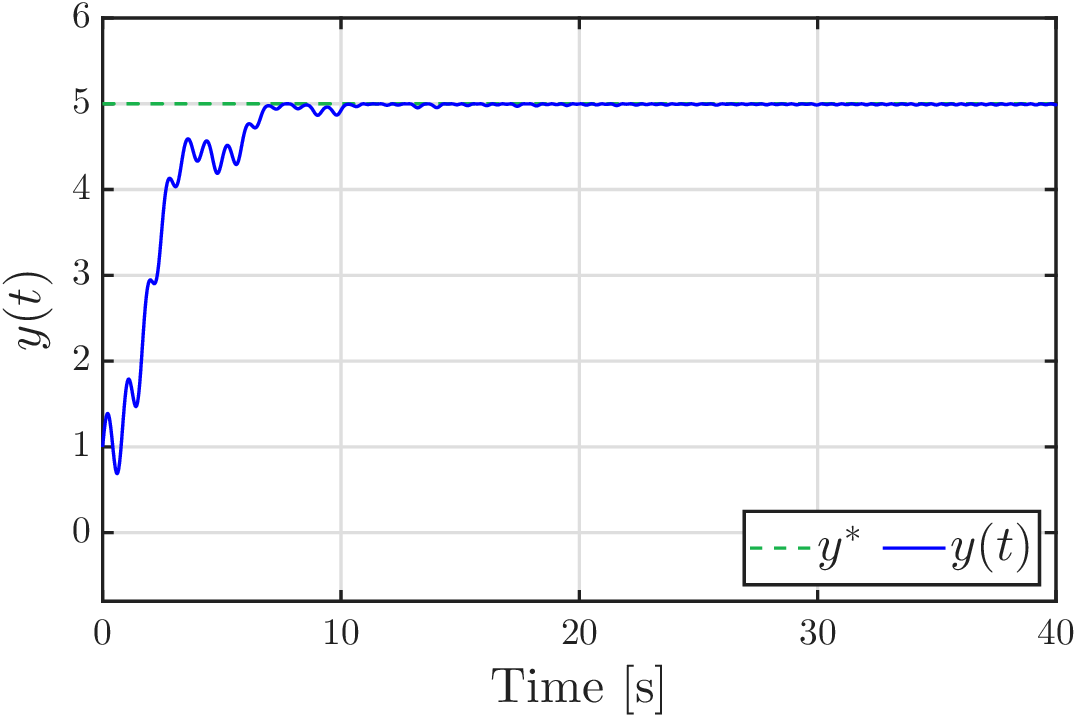}
        \caption{Output $y(t)$ of the static map.}
    \end{subfigure}
    \begin{subfigure}[b]{\columnwidth}
        \centering
        \includegraphics[width=\columnwidth]{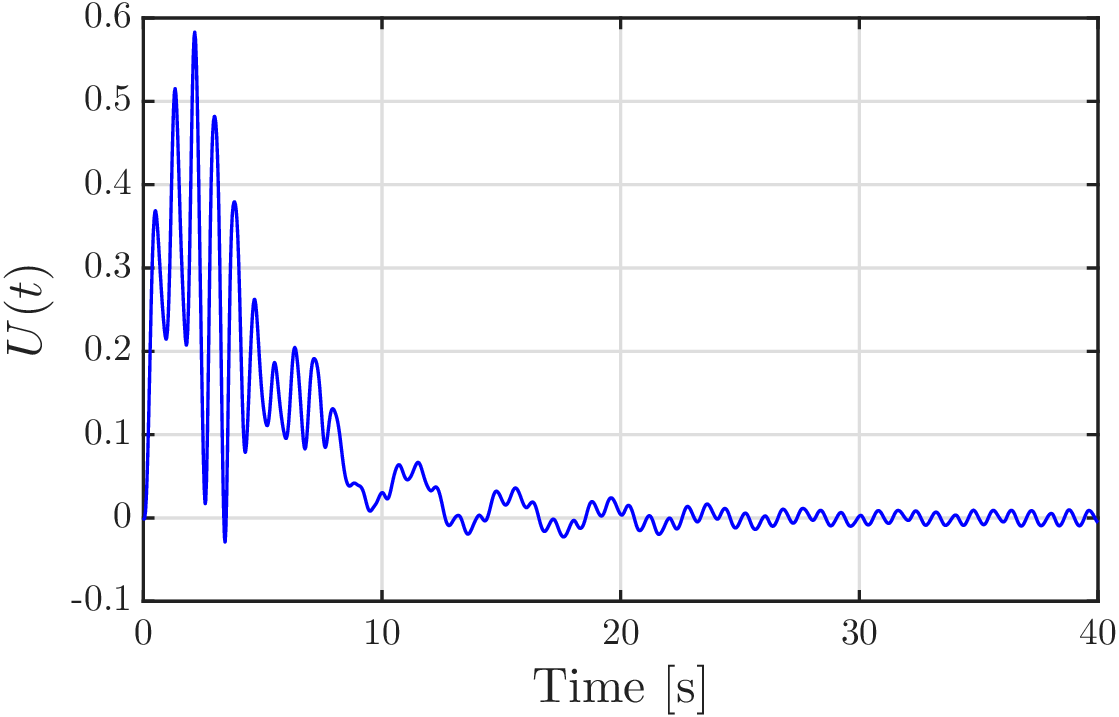}
        \caption{Control signal $U(t)$.}
    \end{subfigure}
    \caption{Closed-loop simulation of the ESC system with actuation dynamics governed by the wave PDE with distributed effect using the controller~\eqref{filteredcontrolreal2026}.}
    \label{fig:closed_loop}
\end{figure}

The output $y(t)$ of the static map  of the ESC system with actuation governed by a wave PDE with distributed effects is shown in Figure~\ref{fig:closed_loop}(a), 
where it is possible to observe the convergence of $y(t)$ to a neighborhood of the unknown extremum point $y^{\ast} = 5$. The signal $U(t)$ of the distributed wave PDE compensation controller implemented as in~\eqref{filteredcontrolreal2026} is presented in Figure~\ref{fig:closed_loop}(b).

\begin{figure}[!ht]
    \centering
    \includegraphics[width=\columnwidth]{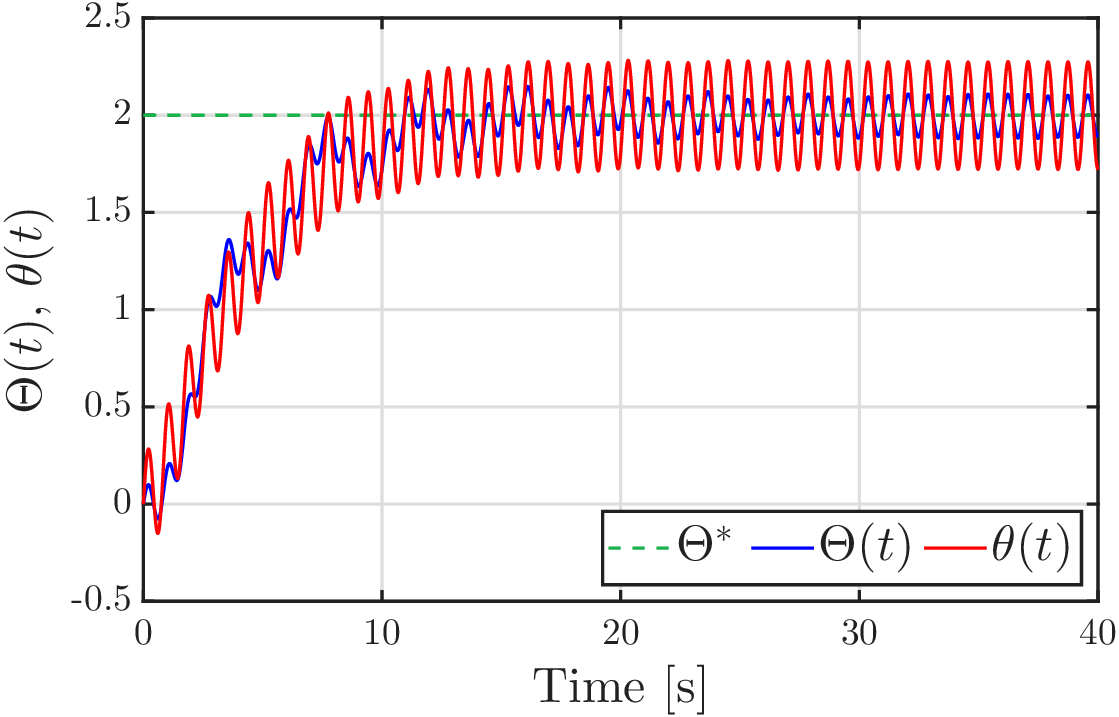}
    \caption{Signals $\theta(t)$ and $\Theta(t)$ of the closed-loop simulation.}
    \label{fig:thetas}
\end{figure}

\newpage
The signals of the actuation input $\theta(t)$ and the distributed actuation input $\Theta(t)$ are shown in Figure~\ref{fig:thetas}. It is also possible to note the convergence of both signals to a neighborhood of the optimizer $\Theta^{\ast} = 2$.

\begin{figure}[!ht]
    \centering
    \begin{subfigure}[b]{\columnwidth}
        \centering
        \includegraphics[width=\columnwidth]{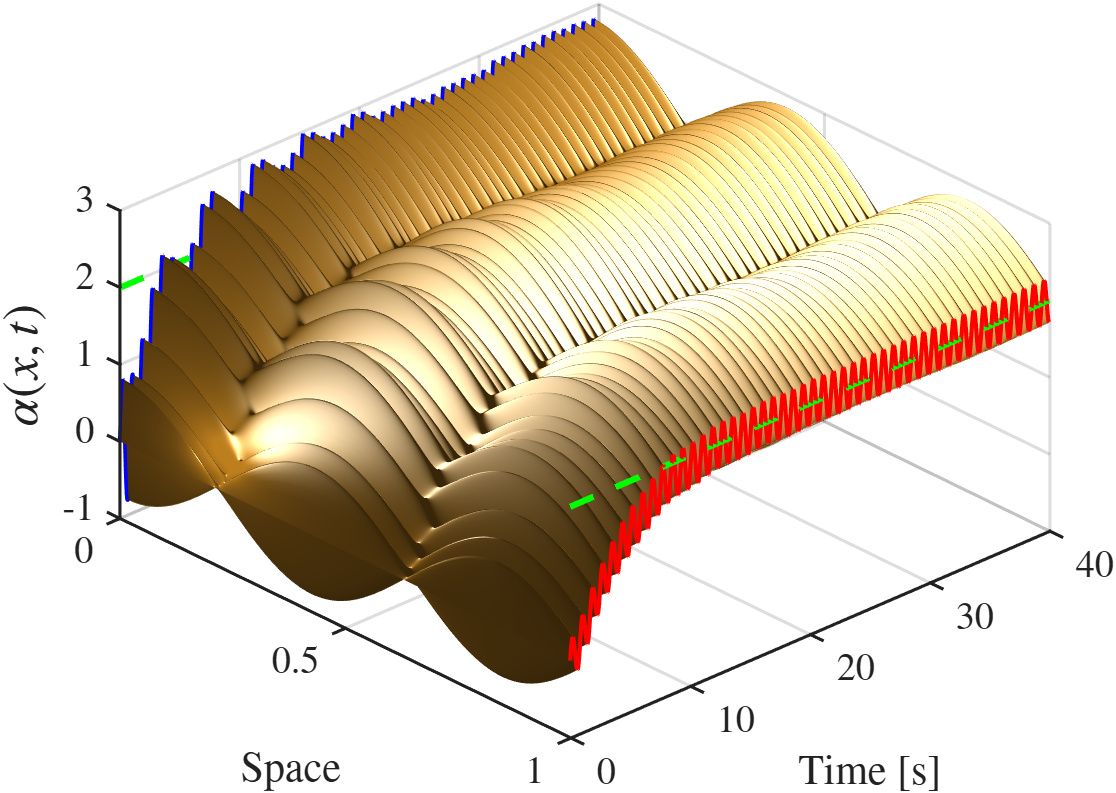}

        \caption{View of $\theta(t) = \alpha(D,t)$ in red.}
    \end{subfigure}
    \begin{subfigure}[b]{\columnwidth}
        \centering
        \includegraphics[width=\columnwidth]{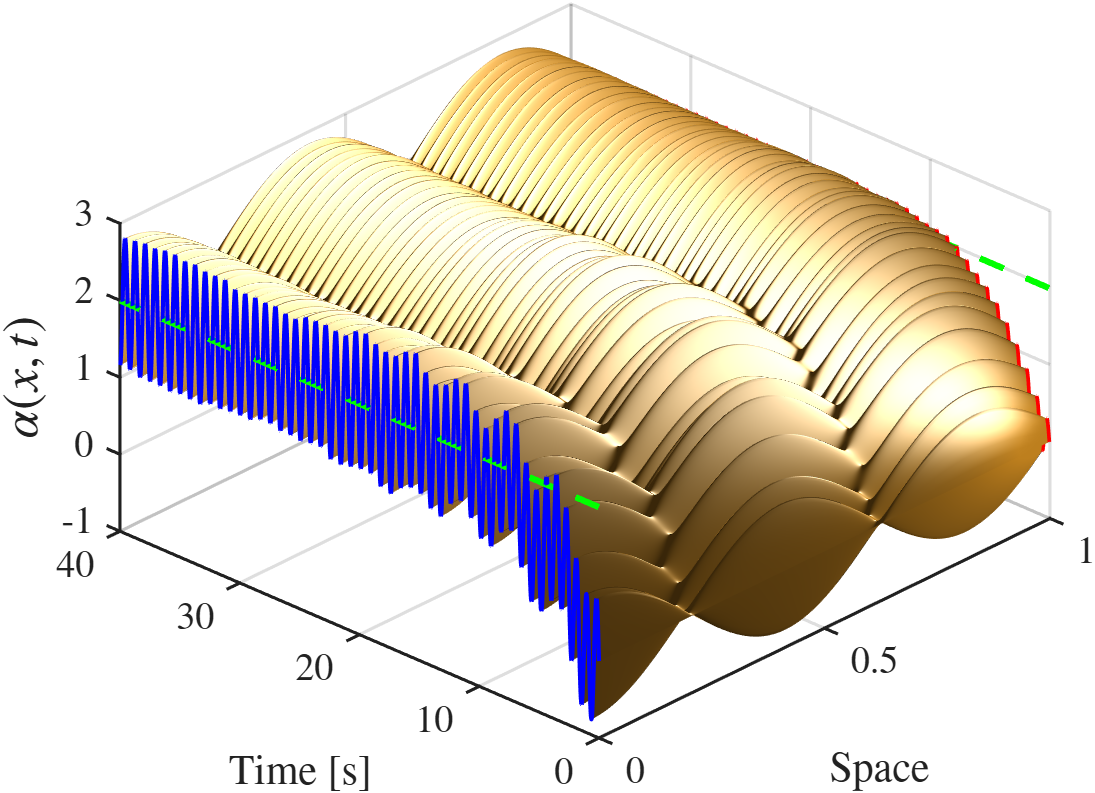}
        \caption{View of $\alpha(0,t)$ in blue.}
    \end{subfigure}
    \caption{Convergence of the state $\alpha(x,t)$ in a three-dimensional space. The signal in red is $\theta(t)=\alpha(D,t)$, while the signal in blue is $\alpha(0,t)$, both reaching a neighborhood of $\Theta^\ast = 2$.}
    \label{fig:3d-view}
\end{figure}

The evolution of the distributed wave PDE~\eqref{eq:6}--\eqref{eq:2} of the closed-loop system in a three-dimensional space with the space domain $x \in [0,1]$ and the time $t$ is shown in Figure~\ref{fig:3d-view}. The curves in blue and red show the convergence of $\alpha(0,t)$ and $\theta(t) = \alpha(D,t)$ 
to a small neighborhood around the optimizer $\Theta^{*} = 2$, respectively. 

\section{Conclusion}
\label{sec:conclusion}
This paper addressed the extremum seeking control problem for distributed wave PDEs with actuation dynamics governed by wave propagation effects. By appropriately defining the estimation errors and designing a distributed wave PDE compensator together with a suitable perturbation signal, the exponential stability of the resulting closed-loop averaged system was rigorously established. Consequently, it was shown that the system trajectories converge to a small neighborhood of the unknown extremum point. Numerical simulations clearly illustrated the effectiveness of the proposed boundary extremm seeking control strategy, even in the presence of this wider class of distributed PDEs.

While related extremum seeking control results exist for systems with distributed delays or diffusion-type dynamics, such approaches cannot be directly extended to the wave PDE setting considered here, due to the fundamentally different hyperbolic nature of the governing equations. In particular, both the control design and the trajectory generation problems addressed in this paper differ substantially from existing methodologies based on delay compensation or parabolic PDE frameworks. The stability analysis developed herein departs from reduction-based techniques and instead relies on a backstepping-inspired design tailored to wave PDE dynamics. Moreover, the proposed trajectory generation mechanism constitutes a novel contribution relative to existing extremum seeking control literature and is specifically suited to handle distributed wave effects. To the best of the authors' knowledge, this work represents one of the first systematic extensions of extremum seeking control to systems with wave PDE actuation with distributed dynamics.

%This paper addressed the problem of extremum seeking control for distributed wave %PDEs. Through the proper definition of the tracking errors, it was possible to %obtain an average system, whose exponential stability was guaranteed by the designed %distributed wave PDE compensator. Thereby, it was demonstrated that the trajectories %converge to a small neighborhood around the optimum point. The numerical simulations %clearly illustrated the effectiveness of the proposed extremum seeking control %strategy, even in the presence of an actuator with dynamics governed by the %distributed wave PDE.

\section*{Acknowledgements}
This work was supported by the Brazilian agencies CNPq (Grant numbers: 309008/2022-0 and 308791/2025-8), CAPES (Finance Code 001), and FAPERJ.

\bibliographystyle{abbrv}

\bibliography{references}

\end{document}